\documentclass[twoside,10pt]{article}
\usepackage{graphicx}
\usepackage{subfigure}
\usepackage{url}
\usepackage{ifpdf}
\usepackage{color}

\ifpdf \DeclareGraphicsExtensions{.pdf,.jpg,.png}
\else  \DeclareGraphicsExtensions{.eps}
\fi

\usepackage{amsfonts}
\usepackage{amscd}
\newcommand\Z{\mathbb Z}
\newcommand\R{\mathbb R}
\newcommand\T{\mathbb T}
\newcommand\ii{{\rm i}\,}
\newcommand\ee{{\rm e}}

\newcommand\eps\varepsilon
\newcommand\A{\mathcal A}

\newcommand\K{\mathcal K}
\newcommand\Lc{\mathcal L}
\newcommand\M{\mathcal M}
\newcommand\Pc{\mathcal P}
\newcommand\Zc{\mathcal Z}

\newcommand\ds{\displaystyle}
\newcommand\dfrac[2]{\ds\frac{#1}{#2}}

\newcommand\p[1]{\left(#1\right)}
\newcommand\pq[1]{\left[#1\right]}
\newcommand\pp[1]{\left\{#1\right\}}
\newcommand\scprod[2]{\left\langle#1,#2\right\rangle}
\newcommand\abs[1]{\left|#1\right|}
\newcommand\absu[1]{\abs{#1}_1}
\newcommand\absd[1]{\abs{#1}_2}

\newcommand\tl{\tilde}
\newcommand\wh[1]{\widehat{#1}}
\newcommand\ol[1]{\overline{#1}}
\newcommand\Ord{{\mathcal O}}
\newcommand\hardsum[2]
  {\ds\sum_{\begin{array}{c}\null^{#1}\\[-3pt]\null^{#2}\end{array}}}

\newcommand\mmatrix[4]{\p{\begin{array}{cc}#1&#2\\[3pt]#3&#4\end{array}}}

\newtheorem{theorem}{Theorem}
\newtheorem{proposition}[theorem]{Proposition}
\newtheorem{lemma}[theorem]{Lemma}

\newcommand\bremark{\noindent{\bf Remark.}\ \ }
\newcommand\eremark{\bigskip}
\newcommand\bremarks{\noindent{\bf Remarks.}\ \bnm}
\newcommand\eremarks{\enm}
\newcommand\proof{\noindent\emph{Proof.}\ \ }
\newcommand\sketchproof{\noindent\emph{Sketch of the proof.}\ \ }
\newcommand\qed{\ \ \null\nolinebreak\hfill$\frame{\large\phantom a}$}

\newcommand\beq{\begin{equation}}
\newcommand\eeq{\end{equation}}
\newcommand\bea{\begin{eqnarray}}
\newcommand\eea{\end{eqnarray}}
\newcommand\bean{\begin{eqnarray*}}
\newcommand\eean{\end{eqnarray*}}
\newcommand\btm{\vspace{-\baselineskip}\begin{itemize}}
\newcommand\etm{\end{itemize}\vspace{-\baselineskip}}
\newcommand\btmm{\begin{itemize}}
\newcommand\etmm{\end{itemize}}
\newcommand\bnm{\vspace{-\baselineskip}\begin{enumerate}}
\newcommand\enm{\end{enumerate}\vspace{-\baselineskip}}

\begin{document}

\title{Exponentially small asymptotic estimates for the splitting of
  separatrices to whiskered tori with quadratic and cubic frequencies
  \ \footnote{This work has been partially supported by the Spanish
      MINECO-FEDER Grants MTM2009-06973, MTM2012-31714
      and the Catalan Grant 2009SGR859.
      The author~MG has also been supported by the
      DFG~Collaborative Research Center TRR~109
      ``Discretization in Geometry and Dynamics''.}}
\author{\sc
    Amadeu Delshams$\,^1$,
  \ Marina Gonchenko$\,^2$,
  \\\sc
  \ Pere Guti\'errez$\,^3$
\\[12pt]
  {\small
  $^{1,3}\;$\parbox[t]{5.2cm}{
    Dep. de Matem\`atica Aplicada I\\
    Universitat Polit\`ecnica de Catalunya\\
    Av. Diagonal 647, 08028 Barcelona\\
    {\footnotesize
      \texttt{amadeu.delshams@upc.edu}\\
      \texttt{pere.gutierrez@upc.edu}}}
  \quad
  $^2\;$\parbox[t]{5.2cm}{
    Tecnische Univesit\"at Berlin\\
    Fakult\"at II -- Mathematik \\und Naturwissenschaften\\
    Institut f\"ur Mathematik, MA 8-18\\
    Stra{\ss}e des 17. Juni 136\\
    D-10623 Berlin\\
    {\footnotesize
      \texttt{gonchenk@math.tu-berlin.de}}}
  }}
\date{May 21, 2003}
\maketitle
\begin{abstract}
We study the splitting of invariant manifolds of whiskered tori with two or
three frequencies in nearly-integrable Hamiltonian systems. We consider
2-dimensional tori with a frequency vector $\omega=(1,\Omega)$ where $\Omega$
is a quadratic irrational number, or 3-dimensional tori with a frequency vector
$\omega=(1,\Omega,\Omega^2)$ where $\Omega$ is a cubic irrational number.
Applying the Poincar\'e--Melnikov method, we find exponentially small
asymptotic estimates for the maximal splitting distance between the stable and
unstable manifolds associated to the invariant torus, showing that such
estimates depend strongly on the arithmetic properties of the frequencies. In
the quadratic case, we use the continued fractions theory to establish a
certain arithmetic property, fulfilled in 24 cases, which allows us to provide
asymptotic estimates in a simple way. In the cubic case, we focus our attention
to the case in which $\Omega$ is the so-called cubic golden number (the real
root of $x^3+x-1=0$), obtaining also asymptotic estimates. We point out the
similitudes and differences between the results obtained for both the quadratic
and cubic cases.
\par\vspace{12pt}
\\[3pt]
\end{abstract}

\section{Introduction}

\subsection{Background and objectives}

The aim of this paper is to introduce a methodology for measuring
the exponentially small
splitting of separatrices in a perturbed Hamiltonian system,
associated to an $\ell$-dimensional whiskered torus
(invariant hyperbolic torus) with an algebraic frequency vector,
quadratic in the case $\ell=2$, and cubic in the case $\ell=3$.

As the unperturbed system, we consider an integrable Hamiltonian $H_0$
with $\ell+1$ degrees of freedom
having $\ell$-dimensional whiskered tori
with coincident stable and unstable whiskers.
In general, for a perturbed Hamiltonian $H=H_0+\mu H_1$ where $\mu$ is small,
the whiskers do not coincide anymore,
giving rise to the phenomenon called \emph{splitting of separatrices},
discovered by Poincar\'e \cite{Poincare90}.
In order to give a measure for the splitting,
one often describes it by a periodic vector function
$\M(\theta)$, $\theta\in\T^\ell$,
usually called \emph{splitting function},
giving the distance between the invariant manifolds
in the complementary directions, on a transverse section.
The most popular tool to measure the splitting is
the \emph{Poincar\'e--Melnikov method},
introduced in \cite{Poincare90} and rediscovered later by
Melnikov and Arnold \cite{Melnikov63,Arnold64}.
This method provides a first order approximation
\beq\label{eq:Melniapprox}
  \M(\theta)=\mu M(\theta)+\Ord(\mu^2),
\eeq
where $M(\theta)$ is called the \emph{Melnikov function}
and is defined by an integral.
In fact, it was established \cite{Eliasson94,DelshamsG00} that both
the splitting and the Melnikov functions are the gradients of
scalar functions: the \emph{splitting potential} and the
\emph{Melnikov potential}, denoted $\Lc(\theta)$ and $L(\theta)$ respectively.
This result implies the existence of \emph{homoclinic orbits}
(i.e.~intersections between the stable and unstable whiskers)
in the perturbed system.

We focus our attention on a concrete torus
with an $\ell$-dimensional frequency vector of \emph{fast frequencies}:
\begin{equation}
  \omega_\varepsilon = \frac\omega{\sqrt{\varepsilon}}\,,
  \quad\eps>0,
\label{eq:omega_eps}
\end{equation}
with a relation between the parameters, of the form $\mu=\eps^p$
for some $p>0$.
Thus, we have a \emph{singular} perturbation problem,
and the interest for this situation lies in its relation
to the normal form in the vicinity
of a simple resonance \cite{Niederman00,DelshamsG01},
of a nearly-integrable Hamiltonian $\K=\K_0+\eps\K_1$.
In such a singular problem, one can give upper bounds for the splitting,
showing that it is \emph{exponentially small} with respect to $\eps$.
The first of such upper bounds was obtained by Neishtadt \cite{Neishtadt84}
in one and a half degrees of freedom i.e.~for 1~frequency,
and later this was extended to the case of~2~or more frequencies
(see for instance
\cite{Simo94,Gallavotti94,BenettinCG97,BenettinCF97,
DelshamsGJS97,DelshamsGS04}).

The problem of establishing \emph{lower bounds}
for the exponentially small splitting,
or even \emph{asymptotic estimates}, is more difficult,
but some results have been obtained also by several methods.
The difficulty lies in the fact that
the Melnikov function is exponentially small in $\eps$ and
the error of the method could overcome the main term
in~(\ref{eq:Melniapprox}).
Then, an additional study is required in order to validate
the Poincar\'e--Melnikov method.
In the case of 1~frequency,
the first result providing an asymptotic estimate for the
exponentially small splitting was obtained
by Lazutin \cite{Lazutkin03} in 1984, for the Chirikov standard map,
using complex parameterizations of the invariant manifolds.
The same technique was used to justify the Poincar\'e--Melnikov method in
a Hamiltonian with one and a half degrees of freedom
\cite{DelshamsS92,DelshamsS97,Gelfreich97}
or an area-preserving map \cite{DelshamsR98}.
In fact, when the Poincar\'e--Melnikov approach cannot be validated,
other techniques can be applied to get exponentially small estimates,
such as complex matching \cite{Baldoma06,OliveSS03,MartinSS11a,MartinSS11b},
or ``beyond all orders'' asymptotic methods \cite{Lombardi00},
or continuous averaging \cite{Treschev97,ProninT00}.

For~2 or more frequencies, it turns out that
\emph{small divisors} appear in the splitting function and,
as first noticed by Lochak \cite{Lochak92},
the arithmetic properties of the frequency vector $\omega$ play an
important r\^ole. This was established by Sim\'o \cite{Simo94},
and rigorously proved in \cite{DelshamsGJS97}
for the quasi-periodically forced pendulum.
A different technique was used by 
Lochak, Marco and Sauzin \cite{Sauzin01,LochakMS03},
and Rudnev and Wiggins \cite{RudnevW00}, namely
the parametrization of the whiskers as solutions of
Hamilton--Jacobi equation, to obtain exponential small estimates
of the splitting, and the existence of transverse homoclinic orbits
for some intervals of the perturbation parameter~$\varepsilon$.
Besides, it was shown in \cite{DelshamsG04} the continuation of
the exponentially small estimates and the
transversality of the splitting, for all sufficiently small values of $\eps$,
under a certain condition on the phases of the perturbation.
Otherwise, homoclinic bifurcations can occur,
studied by Sim\'o and Valls \cite{SimoV01} in the Arnold's example.
The quoted papers considered the case of 2~frequencies,
and assuming in most cases that the frequency ratio
is the famous \emph{golden mean} $\Omega_1=(\sqrt5-1)/2$.
A~generalization to some other quadratic frequency ratios was studied
in \cite{DelshamsG03}.
For a more complete background and references concerning exponentially small
splitting, see for instance \cite{DelshamsGS04}.

The main objective of this paper is to develop a unified methodology
in order to generalize the results on exponentially small splitting
to frequency vectors $\omega$ in~$\R^2$ or~$\R^3$,
in order to obtain asymptotic estimates for the
\emph{maximal splitting distance}
(and, consequently, to show the existence of splitting),
and to emphasize the dependence of such estimates
on the arithmetic properties of the frequencies.
Namely, we consider two possibilities:
\btm
\item \emph{quadratic frequencies}:
  \ $\ell=2$ and $\omega=(1,\Omega)$,
  where $\Omega$ is a quadratic irrational number;
\item \emph{cubic frequencies}:
  \ $\ell=3$ and $\omega=(1,\Omega,\Omega^2)$,
  where $\Omega$ is a cubic irrational number whose two conjugates
  are not real.
\etm
Such frequency vectors satisfy a \emph{Diophantine condition},
\begin{equation}
|\langle k, \omega\rangle| \geq \frac{\gamma}{|k|^{\ell-1}}, \;\;\;
\forall k\in \mathbb{Z}^l\backslash\{0\}
\label{eq:DiophCond}
\end{equation}
with  some $\gamma>0$, in both the quadratic and the cubic cases.
We point out that $\ell-1$ is the minimal possible exponent
for Diophantine inequalities in $\R^\ell$
(see for instance \cite[ap.~4]{LochakM88})

One of the goals of this paper is to show, for the above frequencies,
that we can detect the integer vectors
$k\in\Z^l\setminus\pp0$ providing an approximate equality
in~(\ref{eq:DiophCond}), i.e.~giving the ``least'' small divisors
(relatively to the size of $\abs k$).
We call such vectors $k$ the \emph{primary resonances} of $\omega$,
and other vectors the \emph{secondary} ones.
We show that, if a certain arithmetic condition is fulfilled
(see the separation condition~(\ref{eq:B0A1})),
then the harmonics associated to such primary vectors $k$
are the dominant ones in the splitting function $\M(\theta)$,
for each small enough value of the perturbation parameter $\eps$.

In the quadratic case, the required arithmetic condition~(\ref{eq:B0A1}) can be
formulated in terms of the \emph{continued fraction} of $\Omega$,
which is (eventually) periodic, and
in fact, we can restrict ourselves to the case of \emph{purely periodic}
continued fractions. There are 24 numbers satisfying~(\ref{eq:B0A1}),
all of them having 1-periodic or 2-periodic continued fractions,
\beq\label{eq:24numbers}
  \Omega_a=[\ol a],\ \ a=1,\ldots,13,
  \quad
  \mbox{and}
  \qquad
  \Omega_{1,a}=[\ol{1,a}],\ \ a=2,\ldots,12
\eeq
(this includes the golden number $\Omega_1=[1,1,1,\ldots]=(\sqrt{5}-1)/2$).

In the cubic case, there is no standard continued fraction theory,
but a particular study can be carried out for each cubic irrational $\Omega$.
We consider in this paper the \emph{cubic golden number}
(see for instance \cite{HardcastleK00}):
\beq\label{eq:cubicgolden}
  \Omega\approx0.6823,
  \qquad
  \mbox{the real root of $x^3+x-1=0$,}
\eeq
but we stress that a similar approach could be carried out for other cases.

In the main result of this paper (see Theorem~\ref{thm:main}),
we establish exponentially small asymptotic estimates
for the maximal splitting distance, valid in all
the cases~\mbox{(\ref{eq:24numbers}--\ref{eq:cubicgolden})}.
In this way, we show that the results provided in \cite{DelshamsG03}
for some quadratic frequencies are extended to other cases,
including a particular case of cubic frequencies.
As far as we know, this is the first result
providing asymptotic estimates (and, hence, lower bounds)
for the exponentially small splitting of separatrices with 3~frequencies.
To avoid technicalities, we put emphasis
on the constructive part of the proofs,
using the arithmetic properties of the frequencies
in order to provide a unified methodology which can be applied
to both the quadratic and the cubic cases,
stressing the similarities and differences between them.
We determine, for every $\eps$ small enough,
the dominant harmonic of the Melnikov function $M(\theta)$,
associated to a primary resonance,
and consequently we obtain an estimate for the maximal value of this function.

In a further step, the first order approximation has to be validated
showing that the dominant harmonics of the splitting
function $\M(\theta)$ correspond to the dominant harmonics
of the Melnikov function, as done in \cite{DelshamsG04}.
Besides, one can show in the cases~(\ref{eq:24numbers}--\ref{eq:cubicgolden})
that the invariant manifolds intersect along transverse homoclinic orbits,
with an exponentially small angle.
To obtain this, one needs to consider the ``next'' dominant harmonics
(at least 2~ones in the quadratic case and at least 3~ones in the cubic case,
provided their associated vectors $k$ are linearly independent),
which can be carried out for the frequency vectors considered.
Nevertheless, in some cases the secondary resonances
have to be taken into account giving rise to more involved estimates.
We only provide here the main ideas, and
rigorous proofs will be published elsewhere.

\subsection{Setup and main result}

In order to formulate our main result,
let us describe the Hamiltonian considered,
which is analogous to the one considered
in \cite{DelshamsGS04} and other related works.
In symplectic coordinates
$(x,y,\varphi,I)\in\T\times\R\times\T^l\times\R^\ell$,
\begin{equation}
\begin{array}{l}
H(x,y, \varphi, I) = H_0 (x,y, I) + \mu H_1(x, \varphi),\\[4pt]
H_0 (x, y, I) =
\langle \omega_\varepsilon, I\rangle + \frac{1}{2} \langle\Lambda I, I\rangle
+\dfrac{y^2}2+ \cos x -1,
\quad
H_1 (x, \varphi)= h(x) f(\varphi),
\label{eq:ham1}
\end{array}
\end{equation}
with
\beq
  \label{eq:ham2}
  h(x) = \cos x,
  \qquad
  f(\varphi)= \hardsum{k\in \Zc}{k_l\ge0}
  \ee^{-\rho |k|} \cos(\langle k, \varphi \rangle - \sigma_k),
\eeq
where the restriction in the sum is introduced in order to avoid repetitions.
The Hamiltonian~(\ref{eq:ham1}--\ref{eq:ham2})
is a generalization of the Arnold example 
(introduced in \cite{Arnold64} to illustrate the transition
chain mechanism in Arnold diffusion).
It provides a model for the behavior
of a nearly-integrable Hamiltonian system $\K=\K_0+\eps\K_1$
in the vicinity of a \emph{simple resonance},
after carrying out one step of \emph{resonant normal form}
(see for instance \cite{Niederman00,DelshamsG01}).
In this way, our unperturbed Hamiltonian $H_0$ and the perturbation $\mu H_1$
play the r\^ole of the truncated normal form and the remainder respectively;
in its turn the truncated normal form is an $\Ord(\eps)$-perturbation
of the initial Hamiltonian $\K_0$, making the hyperbolicity appear
(a rescaling leads to the fast frequencies~(\ref{eq:omega_eps})).
The parameters $\eps$ and $\mu$ should not be considered as independent,
but linked by a relation of the type $\mu=\eps^p$.

Notice that the unperturbed system $H_0$
consists of the pendulum given by $P(x,y)= y^2/2 + \cos x -1$ and
$\ell$~rotors with fast frequencies:
$\dot{\varphi}= \omega_\varepsilon + \Lambda I$, $\dot{I}=0$.
The pendulum has a hyperbolic equilibrium at the origin,
and the (upper) separatrix can be parameterized by
$(x_0(s),y_0(s))=(4 \arctan\ee^s,2/\cosh s)$, $s\in\R$.
The rotors system $(\varphi, I)$ has the solutions
$\varphi = \varphi_0+(\omega_\varepsilon + \Lambda I_0)\,t$, $I= I_0$.
Consequently, $H_0$ has an $\ell$-parameter family of $\ell$-dimensional
whiskered invariant tori, with coincident stable and unstable whiskers.
Among the family of whiskered tori,
we will focus our attention on the torus located at $I=0$,
whose frequency vector is $\omega_\varepsilon$ as in~(\ref{eq:omega_eps}),
in our case a quadratic or cubic frequency vector
(for $\ell=2$ or $\ell=3$ respectively).
We also assume the condition of \emph{isoenergetic nondegeneracy}
\beq
    \det \left(
    \begin{array}{cc}
    \Lambda & \omega\\
    \omega^\top & 0
    \end{array}
    \right) \neq 0.
\label{eq:isoenerg}
\eeq
When adding the perturbation $\mu H_1$,
the \emph{hyperbolic KAM theorem} can be applied
(see for instance \cite{Niederman00})
thanks to the Diophantine condition~(\ref{eq:DiophCond})
and the isoenergetic nondegeneracy~(\ref{eq:isoenerg}).
For $\mu$ small enough, the whiskered torus persists
with some shift and deformation, as well as its local whiskers.

In general, for $\mu\ne0$ the (global) whiskers
do not coincide anymore, and one can introduce a \emph{splitting function}
$\M(\theta)$, $\theta\in\T^\ell$,
giving the distance between the whiskers
in the complementary directions, on a transverse section
(\cite[\S5.2]{DelshamsG00}; see also \cite{Eliasson94}).
Applying the Poincar\'e--Melnikov method,
the first order approximation~(\ref{eq:Melniapprox})
is given by the (vector) \emph{Melnikov function} $M(\theta)$.
Both functions $\M(\theta)$ and $M(\theta)$  turn to be gradients
of the (scalar) \emph{splitting potential} $\Lc(\theta)$
and \emph{Melnikov potential} $L(\theta)$, respectively.
The latter one can be defined as follows:
\beq\label{eq:L}
  \ds L(\theta) = - \int_{-\infty}^{\infty}
    [h(x_0(t))-h(0)] f(\theta +\omega_\varepsilon t)\,dt,
  \qquad
  M(\theta) = \nabla L(\theta).
\eeq
Notice that $L(\theta)$ is obtained by integrating $H_1$ along
a trajectory of the unperturbed homoclinic manifold,
starting at the point of the section $s=0$ with phase $\theta$.

In order to emphasize the r\^ole played by the arithmetic properties
of the splitting,
we have chosen for the perturbation the special form given in~(\ref{eq:ham2}).
This form was already considered in \cite{DelshamsG04},
and allows us to deal with the Melnikov function and obtain
asymptotic estimates for the splitting.
Notice that the constant $\rho>0$ in the Fourier expansion
of $f(\varphi)$ in~(\ref{eq:ham2})
gives the complex width of analyticity of this function.
The phases $\sigma_k$ can be chosen arbitrarily for the purpose of this paper.

Now we can formulate our main result, providing asymptotic estimates
for the \emph{maximal splitting distance}
in both the quadratic and cubic cases.

\begin{theorem}[main result]\label{thm:main}
For the Hamiltonian system (\ref{eq:ham1}--\ref{eq:ham2})
with $\ell+1$ degrees of freedom,
satisfying the isoenergetic condition~(\ref{eq:isoenerg}), 
assume that $\varepsilon$ is small enough and $\mu=\varepsilon^p$ with $p>3$.
For $\ell=2$, if $\Omega$ is one of the 24 quadratic
numbers (\ref{eq:24numbers}),
and for $\ell=3$, if $\Omega$ is
the cubic golden number~(\ref{eq:cubicgolden}),
the following asymptotic estimate holds:
$$
\max\limits_{\theta\in \mathbb{T}^\ell} |\M(\theta)|
\sim \frac{\mu}{\varepsilon^{1/\ell}}
\exp \left\{- \frac{C_0 h_1 (\varepsilon)}{\varepsilon^{1/2\ell}}\right\}
$$
where $C_0=C_0(\Omega,\rho)$ is a positive constant, defined in (\ref{eq:C0}).
Concerning the function $h_1(\eps)$,
\btm
\item[\rm(a)]
in the quadratic case $\ell=2$, it is periodic
in $\ln\varepsilon$, with $\min h_1 (\varepsilon) =1$ and
$\max h_1(\varepsilon)= A_1$, with a constant
$A_1=A_1(\Omega)$ defined in~(\ref{eq:quad_A0A1});
\item[\rm(b)]
in the cubic case $\ell=3$, it satisfies the bound
$0< A_0^- \leq h_1 (\varepsilon) \leq A_1^+$,
with constants $A_0^-=A_0^-(\Omega)$ and $A_1^+=A_1^+(\Omega)$
defined in (\ref{eq:cub_A0A1}).
\etm
\end{theorem}

\bremarks
\item
The periodicity in $\ln\eps$ of the function $h_1(\eps)$
in the quadratic case~(a)
was first established in \cite{DelshamsGJS97}
for the quasi-periodically forced pendulum,
assuming that the frequency ratio is the
golden number $\Omega_1$.
Previously, the existence of an oscillatory behavior with
lower and upper bounds had been shown in \cite{Simo94}.
\item
In contrast to the quadratic case, it turns out in the cubic case~(b)
the function $h_1(\varepsilon)$
is not periodic in $\ln\eps$ and has a more complicated  form
(see Figure~\ref{fig:cub_h1}, where one can conjecture
that $h_1(\varepsilon)$ is a quasiperiodic function).
\item
The exponent $p>3$ in the relation $\mu=\eps^p$ can be improved in some
special cases. For instance, if in~(\ref{eq:ham2}) one considers
$h(x)=\cos x-1$, then the asymptotic estimates are valid for $p>2$.
This is related to the fact
that, in this case, the invariant torus remains fixed under the perturbation
and only the whiskers deform \cite{DelshamsGS04}.
\eremarks

This paper is organized as follows. In Section~\ref{sect:arithm}
we study the arithmetic properties of quadratic and cubic frequencies,
and in Section~\ref{sect:asympt_est} we find, for the frequencies
considered in~(\ref{eq:24numbers}--\ref{eq:cubicgolden}),
an asymptotic estimate of the dominant harmonic of the splitting potential,
together with a bound of the remaining harmonics which allows us
to provide an asymptotic estimate for the maximal splitting distance,
as established in Theorem~\ref{thm:main}.

\section{Arithmetic properties of quadratic and cubic frequencies}
\label{sect:arithm}

\subsection{Iteration matrices and resonant sequences}\label{sect:resseq}

We review in this section the technic developed in \cite{DelshamsG03}
for studying the resonances of quadratic frequencies ($\ell=2$),
showing that it admits a direct generalization to
the case of cubic frequencies ($\ell=3$).

In the 2-dimensional case, we consider a \emph{quadratic frequency vector}
$\omega\in\R^2$, i.e.~its frequency ratio is a quadratic irrational number.
Of course, we can assume without loss of generality that the vector has
the form $\omega=(1,\Omega)$.

On the other hand, in the 3-dimensional case we consider
a \emph{cubic frequency vector} $\omega\in\R^3$,
i.e.~the frequency ratios generate a cubic field
(an algebraic number field of degree~3).
In order to simplify our exposition, we assume that the vector
has the form $\omega = (1, \Omega, \Omega^2)$,
where $\Omega$ is a cubic irrational number,
hence the cubic field is $\mathbb{Q}(\Omega)$.

Any quadratic or cubic frequency vector $\omega\in\R^\ell$
satisfies the \emph{Diophantine condition}~(\ref{eq:DiophCond}),
with the minimal exponent~$\ell-1$, see for instance \cite{Cassels57}.
With this in mind, we define the \emph{``numerators''}
\beq\label{eq:defnumerators}
\gamma_k := |\langle k, \omega\rangle |\cdot|k|^{\ell-1},
\qquad
k\in \Zc^l\setminus\{0\},
\eeq
provided a norm $\abs\cdot$ for integer vectors has been chosen
(for quadratic vectors, it was used in \cite{DelshamsG03}
the norm $\abs\cdot=\absu\cdot$,
i.e.~the sum of absolute values of the components of the vector;
but for cubic vectors it will be more convenient to use we use
the Euclidean norm $\abs\cdot=\absd\cdot$).
Our goal is to provide a classification of the integer vectors $k$,
according to the size of $\gamma_k$,
in order to find the primary resonances
(i.e.~the integer vectors $k$ for which $\gamma_k$ is smallest and,
hence, fitting best the Diophantine condition (\ref{eq:DiophCond})),
and study their separation with respect to the secondary resonances.

The key point is to use a result by Koch \cite{Koch99}:
for a vector $\omega\in\R^\ell$ whose frequency ratios generate
an algebraic field of degree~$\ell$,
there exists a \emph{unimodular} matrix $T$
(a square matrix with integer entries and determinant $\pm 1$)
having the eigenvector $\omega$ with associated eigenvalue $\lambda$
of modulus $>1$, and such the other $\ell-1$ eigenvalues are simple
and of modulus $<1$. This result is valid for any dimension~$\ell$,
and is usually applied in the context of renormalization theory
(see for instance \cite{Koch99,Lopesd02a}), since the iteration of the
matrix $T$ provides successive rational approximations to the direction
of the vector~$\omega$.
Notice that the matrix $T$ satisfying the conditions above is not unique
(for instance, any power $T^j$, with $j$ positive, also satisfies them).
We will assume without loss of genericity that $\lambda$ is positive
($\lambda>1$).

In this paper, we are not interested in finding approximations to~$\omega$,
but rather to the \emph{quasi-resonances} of $\omega$,
which lie close to the orthogonal
hyperplane~$\langle\omega\rangle^\bot$.
With this aim, we consider the matrix $U=(T^{-1})^\top$,
which satisfies the following \emph{fundamental equality}:
\begin{equation}
\langle U k, \omega\rangle = \langle k, U^\top \omega\rangle
= \frac{1}{\lambda} \langle k, \omega\rangle.
\label{eq:equalityU}
\end{equation}
We say that an integer vector $k$ is \emph{admissible} if
$|\langle k, \omega\rangle | < 1/2$.
We restrict ourselves to the set $\A$ of admissible vectors, since for any
$k\notin \A$ we have $|\langle k, \omega\rangle | > 1/2$ and
$\gamma_k > |k|^{\ell-1}/2$.
We see from~(\ref{eq:equalityU}) that
if $k\in \A$, then also $U k \in \A$. We say that $k$ is
\emph{primitive} if $k \in \A$ but $U^{-1} k \notin \A$.
We also deduce from~(\ref{eq:equalityU}) that $k$ is primitive if and only if
\beq\label{eq:prim}
  \frac{1}{2\lambda} < |\langle k, \omega\rangle| < \frac{1}{2}.
\eeq

Since the first component of $\omega$ is equal to 1,
it is clear that any admissible vectors can be presented in the form
\[
  \begin{array}{lll}
    k^0(j) = (-\textrm{rint}\,(j \Omega), j),
    &\quad j = \mathbb{Z}\backslash \{0\}
    &\ \mbox{($\ell=2$)},
  \\[4pt]
    k^0(j) = (-\textrm{rint}\,(j_1 \Omega + j_2 \Omega^2), j_1, j_2),
    &\quad j = (j_1,j_2)\in\Z^2\setminus\pp0
    &\ \mbox{($\ell=3$)}
  \end{array}
\]
(we denote $\textrm{rint}\,(x)$ the closest integer to $x$).
If $k^0(j)\in\Z^\ell$ is primitive,
we also say that $j\in\Z^{\ell-1}$ is primitive,
and denote $\Pc$ be the set of such primitives.
Now we define, for each $j\in\Pc$,
the following \emph{resonant sequence} of integer vectors:
\begin{equation}
s(j,n) := U^n k^0(j), \quad n=0,1,2,\ldots
\label{eq:sjn}
\end{equation}
It turns out that such resonant sequences cover the whole set $\A$
of admissible vectors, providing a classification of them.
The properties of such a classification follow
from Proposition~\ref{prop:quadrfreq}
for the case of quadratic frequencies, and from
Proposition~\ref{prop:cubicfreq} for the case of cubic frequencies.

\subsection{Properties of quadratic frequencies}\label{sect:quadr_arithm}

It is well-known that all \emph{quadratic irrational numbers}
$\Omega\in (0,1)$,
i.e. the real roots of quadratic polynomials with rational coefficients,
have the continued fraction
\[
\Omega = [a_1, a_2, a_3, \ldots]\,,
\qquad a_i\in\Z^+,
\]
that is \emph{eventually periodic},
i.e.~periodic starting with some element $a_i$.
In fact, as we see below we can restrict ourselves to the
numbers with \emph{purely periodic} continued fractions
and denote them according to their
periodic part; for an $m$-periodic continued fraction, we write
$\Omega_{a_1,\ldots,a_m}=[\overline{a_1,\ldots,a_m}]$.
For example, the famous \emph{golden number} is
$\Omega_1=[\overline{1}]=({\sqrt{5}-1})/{2}$, and the \emph{silver number}
is $\Omega_2=[\overline{2}]=\sqrt{2}-1$.

For a quadratic frequency $\omega=(1,\Omega)$,
the matrix $T$ provided by Koch's result \cite{Koch99}
can be constructed directly from the continued fraction of $\Omega$.
The quadratic numbers~(\ref{eq:24numbers}), considered in this paper,
have 1-periodic or 2-periodic continued fractions.
Let us write their matrix $T=T(\Omega)$ with $\omega$ as an eigenvector,
and the associated eigenvalue $\lambda=\lambda(\Omega)>1$:
\[
  \begin{array}{lll}
    \mbox{for}\ \Omega=\Omega_a,
    &\qquad T=\mmatrix a110,
    &\quad
     \lambda=\dfrac1{\Omega}\,;
  \\[12pt]
    \mbox{for}\ \Omega=\Omega_{1,a},
    &\qquad T=\mmatrix{a+1}1a1,
    &\quad
     \lambda=\dfrac1{1-\Omega}\,.
  \end{array}
\]

\bremark
In what concerns the contents of this paper, it is enough to consider
quadratic numbers with purely periodic continued fractions,
due to the equivalence of any quadratic number $\wh\Omega$, with an eventually
periodic continued fraction, to some $\Omega$ with a purely periodic one:
$\ds\widehat{\Omega}= \frac{c+ d \Omega}{a+ b \Omega}\,,$
with integers $a$, $b$, $c$, $d$ such that $ad-bc=\pm1$.
Then, it can be shown that the same results apply to both numbers
$\Omega$ and $\wh\Omega$ for $\eps$ small enough.
For instance, the results for the golden number $\Omega_1$
also apply to the \emph{noble} numbers
$\wh\Omega=[b_1,\ldots,b_n,\ol1]$.
We point out that the treshold in $\eps$ of validity of the results,
not considered in this paper, would depend on the non-periodic part
of the continued fraction.
\eremark

Now we consider the resonant sequences defined in~(\ref{eq:sjn}).
For the matrix $T$, let $v_2$ be a second eigenvector
(with eigenvalue $\sigma/\lambda$ of modulus~$<1$,
where $\sigma=\det T=\pm1$);
hence $\omega$, $v_2$ are a basis of eigenvectors.
For the matrix $U=(T^{-1})^\top$, let $u_1$, $u_2$ be a basis of eigenvectors
with eigenvalues $1/\lambda$ and $\sigma\lambda$ respectively.
It is well-known that $\scprod{u_2}\omega=\scprod{u_1}{v_2}=0$.
For any primitive integer $j\in\Pc$, we define the quantities
\[
    r_j:=\scprod{k^0(j)}\omega,
    \qquad
    p_j:=\scprod{k^0(j)}{v_2}.
\]
The properties of the quadratic frequencies
can be summarized in the following
proposition, whose proof is given in~\cite{DelshamsG03}.

\begin{proposition}\label{prop:quadrfreq}
For any primitive $j\in\Pc$, there exists the limit
\[
\gamma^*_j=\lim\limits_{n\to\infty} \gamma_{s(j,n)}
= \abs{r_j} K_j,
\qquad
  K_j=\abs{\frac{p_j}{\scprod{u_2}{v_2}}\,u_2}
  =\abs{k^0(j)-\frac{r_j}{\scprod{u_1}\omega}\,u_1},
\]
and one has:
\btm
\item[\emph{(a)}]\ $ \gamma_{s(j,n)}
= \gamma_j^* + \mathcal{O}(\lambda^{-2 n})$,
\quad $n\geq 0$;
\item[\emph{(b)}]\ $|s(j,n)|
= K_j\,\lambda^n + \mathcal{O}(\lambda^{-n})$,
\quad $n\geq 0$;
\item[\emph{(c)}]\ $\gamma^*_j>\ds\frac{(1+\Omega)|j|-a}{2\lambda}$\,,
\quad
$a= \dfrac{1}{2}\left(1+\dfrac{|u_1|}{|\langle u_1, \omega\rangle|}\right)$.
\etm
\end{proposition}

Since the lower bounds~(c) for the \emph{``limit numerators''} $\gamma^*_j$
are increasing with respect to the primitive $j$,
we can select the minimal of them, corresponding to
some primitive $j_0$. We denote
\begin{equation}
\gamma^*:=\liminf\limits_{|k|\to \infty} \gamma_k
= \min\limits_{j\in \Pc} \gamma^*_j = \gamma^*_{j_0}>0.
\label{ed:quad_gamj0}
\end{equation}
The corresponding sequence $s_0(n):=s(j_0, n)$
gives us \emph{the primary resonances},
and we call \emph{secondary resonances} the integer vectors
belonging to any of the remaining sequences $s(j,n)$, $j\neq j_0$.

We introduce \emph{normalized numerators} $\tl\gamma_k$ and their limits
$\tl\gamma^*_j$, $j\in\Pc$, after dividing by $\gamma^*$,
and in this way $\tl\gamma^*_{j_0}=1$.
We also define a parameter $B_0$ measuring the \emph{separation}
between primary and secondary resonances:
\beq\label{eq:quad_B0}
  \tl\gamma_k:=\frac{\gamma_k}{\gamma^*},
  \qquad
  \tilde{\gamma}^*_{j}:=\dfrac{\gamma^*_j}{\gamma^*},
  \qquad
  B_0:=\min\limits_{j\in\Pc\setminus\{j_0\}}\p{\tilde{\gamma}^*_j}^{1/2},
\eeq
where we included the square root for convenience, see~(\ref{eq:B0A1}).
We are implicitly assuming the hypothesis
that the primitive $j_0$ is unique,
and hence $B_0>1$. In fact, this happens for all the cases we
have explored.

\subsection{Properties of cubic frequencies}\label{sect:cubic_arithm}

Now, we consider a frequency vector of the form
$\omega = (1, \Omega, \Omega^2)$,
where $\Omega$ is a cubic irrational number $\Omega$.
If we consider the matrix $T$ given by Koch's result \cite{Koch99},
mentioned in Section~\ref{sect:resseq},
we can distinguish two possible cases for its three eigenvalues
$\lambda$, $\lambda_2$, $\lambda_3$
(recall that $\lambda>1$ is the eigenvalue with eigenvector $\omega$):
\btm
\item \emph{the real case}: the three eigenvalues
$\lambda$, $\lambda_2$, $\lambda_3$ are real;
\item \emph{the complex case}: only the eigenvalue $\lambda$ is real,
and the other two ones $\lambda_2$, $\lambda_3$
are a pair of complex conjugate numbers.
\etm
These two cases are often called \emph{totally real}
and \emph{non-totally real} respectively.
In this paper we only consider cubic frequency vectors in the complex
or non-totally real case.

\bremark
The reason to restrict ourselves to the complex case is that
the remaining two (complex) eigenvalues have the same modulus.
As we see below, it is natural to extend the results for quadratic frequencies
to cubic frequencies of complex type.
Instead, the study of the real case would require a different approach,
since the behavior of the associated small divisors turns out to be
different from the complex case considered here.
\eremark

Unlike the 2-dimensional quadratic frequencies,
in the case of 3-dimensional cubic frequencies there is no standard
theory of continued fractions
providing a direct construction of the matrix $T$
(however, there are some multidimensional continued fractions algorithms,
which applied to the pair $(\Omega,\Omega^2)$ could be helpful to provide $T$,
see for instance \cite{HardcastleK01,KhaninLM07}).
Fortunately, for a given concrete cubic frequency vector it is not hard to find
the matrix $T$ by inspection,
as we do in Section~\ref{sect:cubicgolden} for the cubic golden number.
Other examples of cubic frequencies and their associated matrices
are given in \cite{Chandre02}
(see also \cite{Lochak92} for an account of examples
and results concerning cubic frequencies).

As in Section~\ref{sect:quadr_arithm}, we are going to establish the
properties of the resonant sequences~(\ref{eq:sjn}).
Let us consider a basis of eigenvectors of $T$, writing the two complex ones
in terms of real and imaginary parts: $\omega$, $v_2+\ii v_3$, $v_2-\ii v_3$,
with eigenvalues $\lambda$, $\lambda_2$ and $\lambda_3=\ol\lambda_2$
respectively. Notice that $\abs{\lambda_2}=\lambda^{-1/2}$;
we denote $\phi:=\arg(\lambda_2)$.

In a similar way, we consider for the matrix
$U = \left(T^{-1}\right)^{\top}$ a basis $u_1$, $u_2+\ii u_3$, $u_2-\ii u_3$
with eigenvalues $\lambda^{-1}$, $\lambda_2^{\,-1}$ and
$\lambda_3^{\,-1}=\ol\lambda_2^{\,-1}$ respectively.
In this way, we avoid working with complex vectors.
One readily sees that
$\scprod{u_2}\omega=\scprod{u_3}\omega=0$,
i.e.~$u_2$ and $u_3$ span the resonant plane $\langle\omega\rangle^\bot$.
Other useful equalities are:
$\scprod{u_1}{v_2}=\scprod{u_1}{v_3}=0$,
$\scprod{u_2}{v_2}=-\scprod{u_3}{v_3}$,
$\scprod{u_2}{v_3}=\scprod{u_3}{v_2}$.
We define $Z_1$, $Z_2$ and $\theta$ through the formulas
\beq\label{eq:Ztheta}
  \frac{1}{2}  (|u_2|^2 + |u_3|^2)=Z_1,
  \qquad
  \frac12(|u_2|^2-|u_3|^2)=Z_2\cos\theta,
  \qquad
  \scprod{u_2}{u_3}=Z_2\sin\theta.
\eeq

For any primitive $j$, we define the quantities
\beq\label{eq:rpq}
    r_j:=\scprod{k^0(j)}\omega,
    \quad
    p_j:=\scprod{k^0(j)}{v_2},
    \quad
    q_j:=\scprod{k^0(j)}{v_3},
\eeq
and $E_j$, $\psi_j$ through the formulas
\beq\label{eq:Epsi}
  \dfrac{\langle v_2, u_2\rangle p_j + \langle v_2, u_3\rangle q_j}
  {\langle v_2, u_2\rangle^2 + \langle v_2, u_3\rangle^2}
  =E_j\cos\psi_j,
  \qquad
  \dfrac{\langle v_2, u_3\rangle p_j
  - \langle v_2, u_2\rangle q_j}
  {\langle v_2, u_2\rangle^2 + \langle v_2, u_3\rangle^2}
  =E_j\sin\psi_j.
\eeq
The following proposition extends the results, given in
Proposition~\ref{prop:quadrfreq} for the quadratic case,
to the complex cubic case.

\begin{proposition}
For any primitive $j=(j_1,j_2)\in\Pc$, the sequence of numerators
$\gamma_{s(j,\cdot)}$ oscillates as $n\to \infty$ between two values,
\beq\label{eq:gamjpm}
  \gamma_j^-=\gamma_j^*\,(1-\delta),
  \qquad
  \gamma_j^+=\gamma_j^*\,(1+\delta),
\eeq
where we define
\[
  \gamma_j^*=\abs{r_j}K_j,
  \qquad
  K_j=E_j^{\,2}Z_1,
  \qquad
  \delta=\frac{Z_2}{Z_1}<1.
\]
We also have:
\btm
\item[\rm(a)]
$\gamma_{s(j,n)} = \gamma^*_j\,(1+\delta\cos [2n\phi + 2\psi_j - \theta])+
\Ord(\lambda^{-3n/2})$;
\item[\rm(b)]
$|s(j,n)|^2 =K_j\,(1+\delta\cos [2n\phi + 2\psi_j - \theta])\cdot\lambda^n
  +\Ord(\lambda^{-n/2})$;
\item[\rm(c)]
$\ds\gamma_j^- \geq \frac{1-\delta}{2\lambda(1+\delta)}
\left[ |j| - \frac{|u_1|}{2|\langle u_1, \omega\rangle|}\right]^2$.
\etm
\label{prop:cubicfreq}
\end{proposition}

\proof
We present the primitive vector associated to $j$ in the basis
$u_1$, $u_2$, $u_3$:
$$
k^0(j) = c_1 u_1 + c_2 u_2 + c_3 u_3,
$$
and taking scalar products with $\omega$, $v_2$ and $v_3$
and solving a linear system,
one can obtain the values of the coefficients:
\begin{equation}
c_1 = \frac{r_j}{\langle u_1,\omega\rangle},
\qquad
c_2=E_j\cos\psi_j,
\qquad
c_3=E_j\sin\psi_j,
\label{eq:cj1}
\end{equation}
where the definitions~(\ref{eq:rpq}--\ref{eq:Epsi})
have been taken into account.
Now, we apply the iteration matrix $U$.
Using the identities
\[
  U^nu_2=\lambda^{n/2}[\cos(n\phi)\,u_2+\sin(n\phi)\,u_3],
  \quad
  U^nu_3=\lambda^{n/2}[-\sin(n\phi)\,u_2+\cos(n\phi)\,u_3],
\]
we find
\[
  s(j,n)
  =U^n k^0(j)
  =\lambda^{n/2}E_j[\cos(n\phi+\psi_j)\,u_2+\sin(n\phi+\psi_j)\,u_3]
   +\Ord(\lambda^{-n}),
\]
and we deduce, according to the definitions~(\ref{eq:Ztheta}),
$$
|s(j,n)|^2 =  \lambda^n E_j^{\,2} (Z_1 + Z_2 \cos [2n\phi
+ 2\psi_j - \theta]) + \mathcal{O}(\lambda^{-n/2}),
$$
which gives~(b).
Multiplying by
$|\langle s(j,n), \omega\rangle | = \abs{r_j}\lambda^{-n}$,
we obtain $\gamma_{s(j,n)}$ as given in~(a).
This implies the asymptotic bounds introduced in~(\ref{eq:gamjpm}).

%----------
\begin{figure}[!b]
  \centering
  \includegraphics[width=0.7\textwidth]{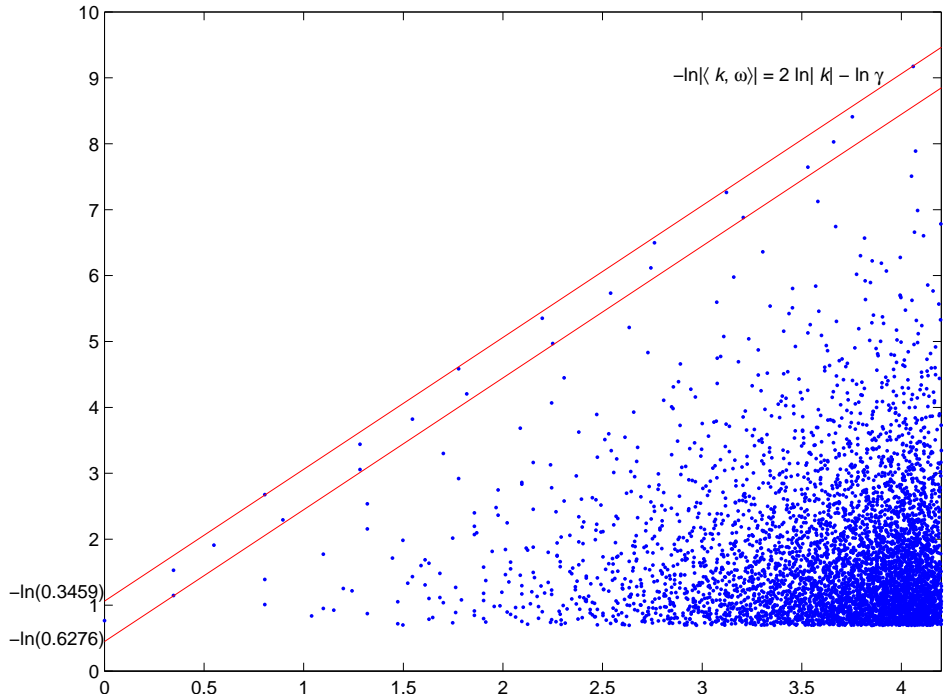}
  \caption{{\small\emph{Points $(\ln|k|, -\ln|\langle k, \omega\rangle|)$.}}}
\label{fig:ln_gam}
\end{figure}
%----------

Finally, one easily sees that
\[
  |k^0(j) - c_1 u_1|^2
  = |c_2 u_2 + c_3 u_3|^2
  =E_j^{\,2}(Z_1+Z_2\cos[2\psi_j-\theta])
  \le E_j^{\,2}(Z_1+Z_2),
\]
and hence using (\ref{eq:cj1}) and (\ref{eq:prim}),
and also that $\abs{k^0(j)}\ge\abs j$, we get
\[
  K_j=E_j^{\,2}Z_1
  \ge \frac{|k^0(j) - c_1 u_1|^2}{1+\delta}
  \ge \frac{1}{1+\delta} \left[ |j| -
  \frac{|u_1|}{2|\langle u_1, \omega\rangle|}\right]^2,
\]
which implies the lower bound given in~(c).
\qed

As we can see in~(a), the existence of limit of the
sequences $\gamma_{s(j,n)}$ stated in Proposition~\ref{prop:quadrfreq}
for the quadratic case, is replaced here by an
\emph{oscillatory limit behavior}, with a lower limit
$\ds\liminf_{n\to\infty}\gamma_{s(j,n)}\ge\gamma_j^-$
and an upper limit $\ds\limsup_{n\to\infty}\gamma_{s(j,n)}\le\gamma_j^+$.
Notice that such values of the limits are exact if the phase $\phi/2\pi$,
that appears in~(a), is irrational.

Selecting the primitive $j_0\in\Pc$ which gives the minimal limits,
we have the \emph{primary resonances},
and we denote them by $s_0(n):= s(j_0,n)$,
and we call \emph{secondary resonances} the integer vectors
belonging to any of the remaining sequences $s(j,n)$, $j\neq j_0$.
Such primary resonances can easily
be detected thanks to Proposition \ref{prop:cubicfreq}(c): although
$\gamma^\pm_j$ are not increasing in general with respect to $\abs j$,
we have an increasing lower bound, which implies that
$\lim\limits_{|j|\to \infty} \gamma^\pm_j = \infty$, and then one
has to check only a finite number of primitive vectors $j$ in order
to find the minimal $\gamma_j^-$ and $\gamma_j^+$ and, hence,
the primary resonances.

As in Section~\ref{sect:quadr_arithm}, we define normalized values
$\tl\gamma_k$, $\tl\gamma_j^*$, $\tl\gamma_j^\pm$, after dividing by
the minimal among the values $\gamma_j^*$,
\begin{equation}\label{eq:cub_gamj0}
  \tl\gamma_k:=\frac{\gamma_k}{\gamma^*},
  \qquad
  \tl\gamma_j^*:=\frac{\gamma_j^*}{\gamma^*},
  \qquad
  \tl\gamma_j^\pm:=\frac{\gamma_j^\pm}{\gamma^*},
  \qquad
  \textrm{where
    \ $\gamma^*:= \min\limits_{j\in \Pc} \gamma^*_j = \gamma^*_{j_0}$.}
\end{equation}
We also introduce a parameter $B_0^-$,
as a measure for the \emph{separation} between primary and secondary
resonances:
\beq\label{eq:cub_B0}
  B_0^-:=\min_{j\in\Pc\setminus\pp{j_0}}
    \p{\frac{\tl\gamma^-_j}{\tl\gamma^+_{j_0}}}^{1/3}
\eeq
(compare with~(\ref{eq:quad_B0}), and see also~(\ref{eq:B0A1})).
Notice that the distinction between primary and secondary resonances
makes sense if $B_0^->1$,
i.e.~the interval $[\tl\gamma^-_{j_0},\tl\gamma^+_{j_0}]$
has no intersection with any other interval
$[\tl\gamma^-_j,\tl\gamma^+_j]$, $j\ne j_0$
(as happens in the cubic golden case, see the next section).

%----------
\begin{table}[!b]
\centering
\begin{tabular}{|c||c|c|c|}
\hline
$k_0(j)$& $\gamma_j^-$ & $\gamma_j^*$& $\gamma_j^+$\\
\hline \hline
$[0, 0, 1]$  &0.3459      &0.4867 &0.6276 \\
$[-1, 2, 0]$ &1.0376      &1.4602 &1.8829 \\
$[-2, 1, 2]$ &3.1127      &4.3807 &5.6488 \\[4pt]
$|j|\ge3$    &$\ge1.2742$ &       &       \\
\hline
\end{tabular}
\caption{{\small\emph{Numerical data for the cubic golden frequency vector}}}
\label{table:boundsgam1}
\end{table}
%----------

\subsection{The cubic golden frequency vector}\label{sect:cubicgolden}

Now, we assume that $\Omega$ is the \emph{cubic golden number}:
the real root of $x^3 + x - 1=0$. We have $\Omega \approx 0.6823$.
In this case, the matrix $T$ can easily be found by inspection.
We have
$$
T=\left(
\begin{array}{ccc}
1& 0& 1\\
1& 0& 0\\
0& 1& 0\\
\end{array}
\right),
\qquad
U= (T^{-1})^\top=\left(
\begin{array}{ccc}
0& 0& 1\\
1& 0& -1\\
0& 1& 0\\
\end{array}
\right),
$$
with the eigenvalue $\lambda=1/\Omega\approx1.4656$.

It is not hard to compute the data provided by Proposition~\ref{prop:cubicfreq}
in this concrete case.
In particular, we have
\beq\label{eq:phi}
  \phi=\arg(\lambda_2)
  =-\arctan \dfrac{4\sqrt{31}}{\Omega(6\Omega^2+9\Omega+4)}+\pi
  \ \approx\ \dfrac{13 \pi}{22}
\eeq
and, from Proposition~\ref{prop:cubicfreq}(b), we have the following
approximately periodic behaviors:
$\gamma_{s(j,n+22)}\approx\gamma_{s(j,n)}$,
and $\abs{s(j,n+22)}\approx\lambda^{11}\abs{s(j,n)}$.
Other relevant parameters are:
$\gamma^*=\dfrac2{31}(5+\Omega+4\Omega^2)\approx0.4867$
and
$\delta=(3-2 \Omega)\sqrt{2-\Omega+\Omega^2}/(5+\Omega+4\Omega^2)
  \approx0.2895$.
In Table~\ref{table:boundsgam1}, we write down the values $\gamma^*_j$,
as well as the bounds $\gamma^-_j$ and $\gamma^+_j$,
for the resonant sequences induced by a few primitives $k_0(j)$,
and a lower bound for all other primitives.
The smallest ones correspond to the primitive
vector $k_0([0, 1])= [0, 0, 1]$ (primary resonances).
The parameter introduced in~(\ref{eq:cub_B0}),
indicating the separation between the primary
and the secondary resonances, is $B_0^-\approx1.1824$.

Additionally, it is interesting to visualize
such a separation in the following way.
Taking logarithm in both hands of the Diophantine
condition (\ref{eq:DiophCond}), we can write it as
$$-\ln |\langle k, \omega \rangle| \le 2 \ln |k| - \ln \gamma.$$
If we draw all the points with coordinates
$(\ln |k|, -\ln |\langle k, \omega \rangle|)$
(see Figure \ref{fig:ln_gam}), we can see a sequence of
points lying between the two straight lines
$-\ln |\langle k, \omega \rangle| = 2 \ln |k| - \ln \gamma^\pm_{[0,1]}$.
Such points correspond to integer vectors belonging
to the sequence of primary resonances: $k=s_0(n)$, $n\ge0$.

\section{Asymptotic estimates for the maximal splitting distance}
\label{sect:asympt_est}

In order to provide asymptotic estimates (or lower bounds) for the splitting,
we start with the first order approximation, given by the
Poincar\'e--Melnikov method.
It is convenient for us to work with the (scalar) Melnikov potential $L$
and the splitting potential $\Lc$, but we state our main result in terms
of the splitting function $\mathcal{M}=\nabla\mathcal{L}$,
which gives a measure of the splitting distance
between the invariant manifolds of the whiskered torus.
Notice also that the nondegenerate critical points of $\Lc$ correspond
to simple zeros of $\mathcal{M}$, and give rise
to transverse homoclinic orbits to the whiskered torus.

In the present paper, we restrict ourselves
to present the constructive part of the proofs,
which corresponds to find, for every sufficiently small $\eps$,
the dominant harmonics of the Fourier expansion
of the Melnikov potential $L(\theta)$,
as well as to provide bounds for the sum
of the remaining terms of that expansion.
The final step, to ensure that the Poincar\'e--Melnikov
method~(\ref{eq:Melniapprox})
predicts correctly the size of splitting
in our singular case $\mu=\varepsilon^p$,
can be worked out simply by
showing that the asymptotic estimates of the dominant harmonics
are large enough to overcome the harmonics of the error term.
This final step is analogous to the one done
in \cite{DelshamsG04} for the case of the golden number~$\Omega_1$
(using the upper bounds for the error term provided in \cite{DelshamsGS04}),
and will be published elsewhere for all cases considered
in Theorem~\ref{thm:main}.

First, we are going to find in Section~\ref{sect:dominant}
an exponentally small asymptotic estimate for the dominant harmonic
among the ones associated to primary resonances,
given by a function $h_1(\eps)$ in the exponent.
We also provide an estimate for the sum of all other
(primary or secondary) harmonics.
This can be done jointly for both
the quadratic ($\ell=2$) and cubic ($\ell=3$) cases.
In Section~\ref{sect:h1}, we establish
a condition ensuring that the dominant harmonic among all harmonics
is given by a primary resonance. This condition is fulfilled
for the frequencies~(\ref{eq:24numbers}--\ref{eq:cubicgolden}).
To complete the proof of Theorem~\ref{thm:main},
we show that the different arithmetic properties
of quadratic and cubic frequencies lead to different properties of
the function $h_1(\eps)$:
periodic (with respect to $\ln\eps$) in the quadratic case,
and a more complicated bounded function in the cubic case.

\subsection{Dominant harmonics of the splitting potential}\label{sect:dominant}

We put our functions $f$ and $h$, defined in~(\ref{eq:ham2}),
into the integral (\ref{eq:L}) and get the Fourier expansion of the Melnikov
potential,
\[
L(\theta) = \sum\limits_{k\in \mathcal{Z}\backslash \{0\}}
L_k \cos(\langle k, \theta\rangle -\sigma_k),
\qquad
L_k = \frac{2\pi |\langle k, \omega_\varepsilon\rangle|
\,\ee^{-\rho |k|}}{\sinh |\frac{\pi}{2}
\langle k, \omega_\varepsilon\rangle|}\,.
\]
Using~(\ref{eq:omega_eps}) and~(\ref{eq:defnumerators}),
we present the coefficients in the form
\begin{equation}
\label{eq:alphabeta}
L_k = \alpha_k\,\ee^{- \beta_k},
\qquad
\alpha_k \approx \frac{4 \pi\gamma_k}{|k|^{\ell-1}\sqrt{\varepsilon}}\,,
\quad
\beta_k =\rho |k| + \frac{\pi \gamma_k}{2 |k|^{\ell-1}\sqrt{\varepsilon}}\,,
\end{equation}
where an exponentially small term has been neglected in the denominator
of $\alpha_k$.
For any given $\eps$, the harmonics with largest coefficients $L_k(\eps)$
correspond essentially to the smallest exponents
$\beta_k(\eps)$. Thus, we have to study the dependence on $\eps$ of
such exponents.

With this aim, we introduce for any $X$, $Y$ (and a fixed $\ell=2,3$)
the function
\beq\label{eq:defG}
  G(\eps;X,Y):=
  \frac{Y^{1/\ell}}\ell
  \pq{(\ell-1)\p{\frac\eps X}^{1/2\ell}+\p{\frac X\eps}^{(\ell-1)/2\ell}}.
\eeq
One easily checks that this function has its minimum at
$\varepsilon=X$, and the corresponding minimum value is $G(X;X,Y)=Y^{1/\ell}$.
Then, the exponents $\beta_k(\eps)$ in~(\ref{eq:alphabeta}) can be
presented in the form
\[
  \beta_k(\eps)=\frac{C_0}{\eps^{1/2\ell}}\,g_k (\varepsilon),
\]
where we define
\bea
  \label{eq:gk}
  &&g_k(\eps):=G(\eps;\eps_k,\tl\gamma_k),
  \qquad
  \eps_k:= D_0\,\frac{\tl\gamma_k^{\,2}}{\abs k^{2\ell}}\,,
\\
  \label{eq:C0}
  &&C_0=\ell\p{\frac\rho{\ell-1}}^{(\ell-1)/\ell}
      \cdot\p{\frac{\pi\gamma^*}2}^{1/\ell},
  \qquad
  D_0=\p{\frac{(\ell-1)\pi\gamma^*}{2\rho}}^2,
\eea
with $\gamma^*=\gamma^*_{j_0}$ and $\tl\gamma_k$
given in~(\ref{ed:quad_gamj0}--\ref{eq:quad_B0})
and~(\ref{eq:cub_gamj0}), for $\ell=2$ and $\ell=3$ respectively.
Consequently, for all $k$ we have
$\ds\beta_k \geq \frac{C_0\tl\gamma_k^{1/\ell}}{\varepsilon^{1/2\ell}}$,
which provides the maximum value of the coefficient $L_k(\eps)$
of the harmonic given by the integer vector~$k$.
Recall that for $k=s(j,n)$, belonging to the resonant sequence
generated by a given primitive $j\in\Pc$ (see definition~(\ref{eq:sjn})),
the (normalized) numerators $\tl\gamma_k$
tend to a limit $\tl\gamma_j^*$
(see Proposition~\ref{prop:quadrfreq} for the quadratic case),
or oscillate between two limit values $\tl\gamma_j^\pm$,
with $\tl\gamma_j^-<\tl\gamma_j^*<\tl\gamma_j^+$ 
(see Proposition~\ref{prop:cubicfreq} for the cubic case).

The primary integer vectors $k$, belonging to the sequence $s_0(n)=s(j_0,n)$,
play an important r\^{o}le here, since they give
the smallest limit $\tl\gamma^*_{j_0}=1$.
Consequently, they give the dominant harmonics of the Melnikov potential,
at least for $\eps$ close to their minimum points $\eps_k$.
Our aim is to show that this happens also for any $\eps$ (small enough)
not necessarily close to $\eps_k$,
under a condition ensuring that the separation between the primary
and secondary resonances is large enough
(see the separation condition~(\ref{eq:B0A1})).
For the sequence of primary resonances, the asymptotic behavior
of the functions $g_k(\eps)$ in~(\ref{eq:gk}), $k=s_0(n)$,
is obtained from the main terms, as $n\to\infty$, given by
Propositions~\ref{prop:quadrfreq} and~\ref{prop:cubicfreq}.
Thus, we can write
\beq\label{eq:gn}
  g_{s_0(n)}(\eps)\ \approx\ g^*_n(\eps):=G(\eps;\eps^*_n,b_n),
  \qquad
  \varepsilon^*_n:= \frac{D_0}{K_{j_0}^{\,2\ell/(\ell-1)}}
   \p{\frac{1}{b_n\lambda^{\ell\,n}}}^{2/(\ell-1)},
\eeq
where we define, in order to unify the notation,
\beq\label{eq:bn}
  \begin{array}{ll}
    b_n=1
    &\textrm{in the quadratic case;}
  \\[3pt]
    b_n=1+\delta\cos[2n\phi+2\psi_{j_0}-\theta]
    &\textrm{in the cubic case.}
  \end{array}
\eeq
Notice that~(\ref{eq:gn}) relies in the approximations
$\tl\gamma_{s_0(n)}\approx b_n$ and
$\abs{s_0(n)}^{\ell-1}\approx K_{j_0}b_n\lambda^n$,
as well as the fact that $\tl\gamma^*_{j_0}=1$.

\begin{figure}[!b]
  \centering
  \includegraphics[width=0.7\textwidth]{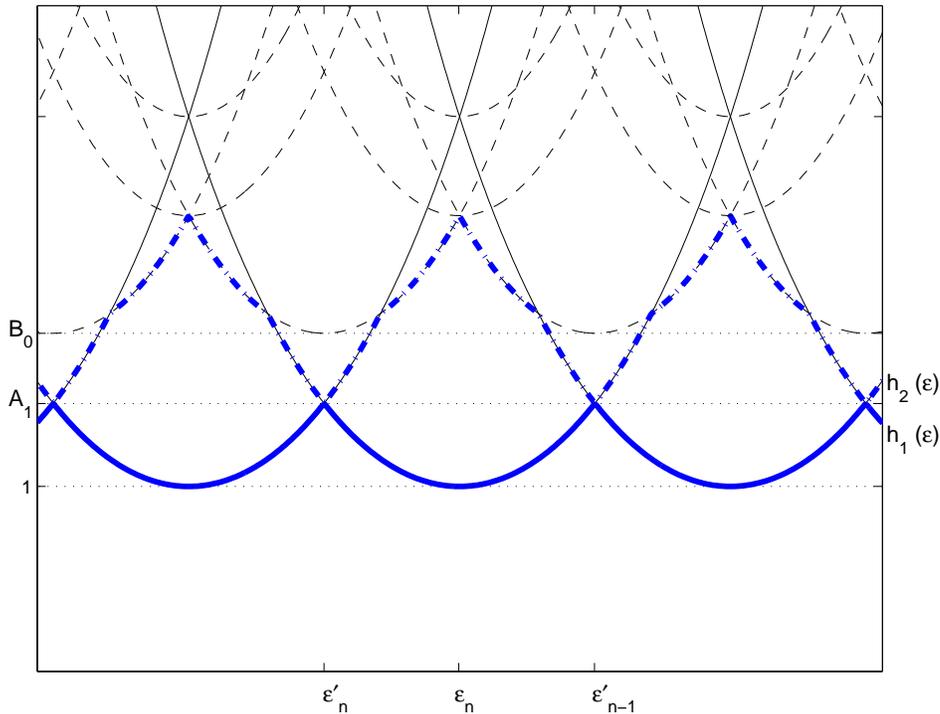}
  \caption{{\small\emph{Graphs of the functions $g_k(\varepsilon)$,
    $h_1(\varepsilon)$, $h_2(\varepsilon)$ for the
    quadratic number $\Omega_{1,2}$,
    using a logarithmic scale for $\varepsilon$}}}
\label{fig:qn-1-2}
\end{figure}

\begin{figure}[!b]
  \centering
  \includegraphics[width=0.7\textwidth]{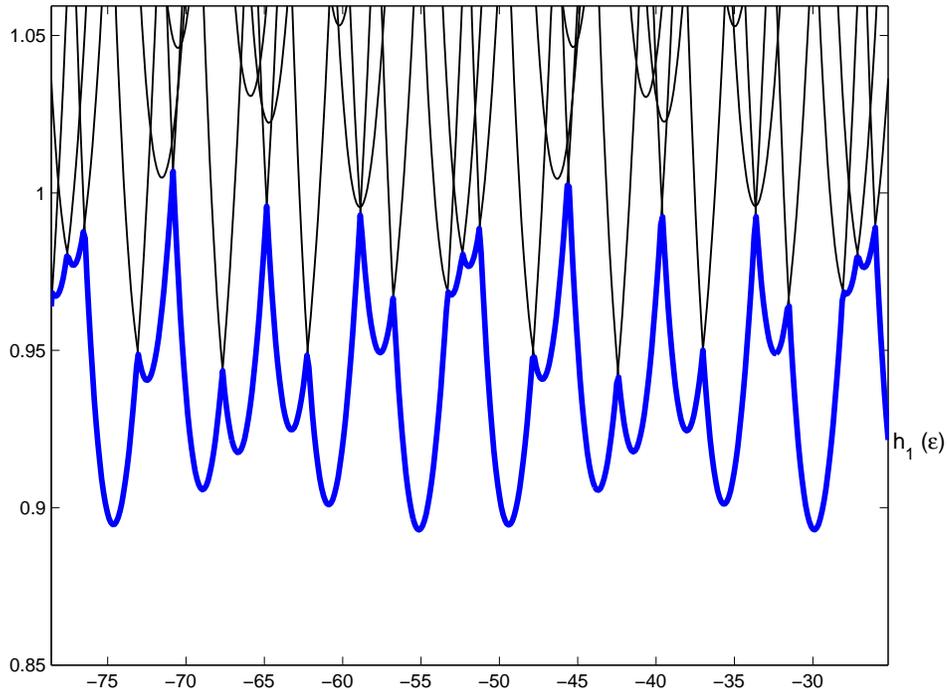}
  \caption{{\small\emph{Graphs of the functions $g^*_n(\varepsilon)$
    and $h_1 (\varepsilon)$ for the cubic golden number $\Omega\approx0.6823$,
    using a logarithmic scale for $\varepsilon$.}}}
\label{fig:cub_h1}
\end{figure}

Now we define, for any given $\eps$, the function $h_1(\eps)$
as the minimum of the values $g^*_n(\eps)$, which takes place
for some index $N=N(\eps)$,
\beq\label{eq:h1}
  h_1(\eps):=\min_{n\ge0}g^*_n(\eps)=g^*_N(\eps).
\eeq
As we show in the next result,
the function $h_1(\eps)$ indicates for any $\eps$
the size of the dominant harmonic among the primary resonances,
given by the integer vector $k=s_0(N)$.
Besides, we are going to establish an
asymptotic estimate for the sum of
all remaning coefficients in the Fourier expansion of the splitting function.
This second estimate is written in terms of the function
\beq\label{eq:h2}
  h_2(\eps):=\min_{k\ne s_0(N)}g_k(\eps).
\eeq
See Figure~\ref{fig:qn-1-2} as an illustration for the functions
$h_1(\eps)$ and $h_2(\eps)$,
corresponding to the case of the quadratic vector
given by $\Omega_{1,2}$, and Figure~\ref{fig:cub_h1}
for $h_1(\eps)$ in the case of the cubic golden vector.

%----------
%----------

Notice that this result is stated in terms of the Fourier coefficients of
the splitting function $\M=\nabla\Lc$.
We write, for the splitting potential,
\[
  \Lc(\theta) = \sum\limits_{k\in \mathcal{Z}\backslash \{0\}}
  \Lc_k \cos(\langle k, \theta\rangle -\tau_k),
\]
with upper bounds for $\abs{\Lc_k-\mu L_k}$ and $\abs{\tau_k-\sigma_k}$.
Then, for the splitting function we have $\abs{\M_k}=\abs k\Lc_k$.

\begin{proposition}\label{prop:Mdom}
For $\eps$ small enough and $\mu=\varepsilon^p$ with $p>3$, one has:
\btm
\item[\rm(a)]
$\ds\abs{\M_{s_0(N)}}
 \sim\mu\abs{s_0(N)}L_{s_0(N)}
 \sim\frac{\mu}{\varepsilon^{1/\ell}}
   \exp\pp{-\frac{C_0h_1(\eps)}{\eps^{1/2\ell}}}$;
\item[\rm(b)]
$\ds\sum_{k\ne s_0(N)}\abs{\M_k}
 \sim\frac{\mu}{\varepsilon^{1/\ell}}
   \exp\pp{-\frac{C_0h_2(\eps)}{\eps^{1/2\ell}}}$.
\etm
\end{proposition}

\sketchproof
At first order in $\mu$, for the coefficients of the splitting function
we can write
\beq\label{eq:Mk}
  \abs{\M_k}\sim\mu\abs kL_k=\mu\abs k\alpha_k\,\ee^{- \beta_k},
\eeq
where we have neglected the error term in the Melnikov
approximation~(\ref{eq:Melniapprox}), and we have used
the expression~(\ref{eq:alphabeta}) for the coefficients of
the Melnikov potential.
As mentioned throughout this section,
the main behavior of the coefficients~$L_k$ is given by
the exponents $\beta_k$, which have been written in~(\ref{eq:gk})
in terms of the functions $g_k(\eps)$.
In particular, the coefficient associated to the dominant harmonic, among
the primary resonances, $L_{s_0(N)}$ with $N=N(\eps)$,
can be expressed in terms of the function $h_1(\eps)$
introduced in~(\ref{eq:h1}).

Now, we consider the remaining factors in~(\ref{eq:Mk}).
We see from~(\ref{eq:alphabeta}) that such factors can be written as
$\abs k\alpha_k\sim\abs k^{-(\ell-2)}\sqrt\eps$.
For $k=s_0(N)$, let us show that they
turn out to be polynomial with respect to $\eps$,
with a concrete exponent to be determined.
First, we use that in~(\ref{eq:gn}) we have
$\eps^*_n\sim\lambda^{-2\ell n/(\ell-1)}$.
Then, for a given $\eps$ the coefficient $N=N(\eps)$ giving the dominant 
harmonic is such that $\eps^*_N$ (the minimum of the function $g^*_N$)
is close to $\eps$, and hence $\lambda^N\sim\eps^{-(\ell-1)/2\ell}$.
On the other hand, we deduce from
from Propositions~\ref{prop:quadrfreq}(b) and~\ref{prop:cubicfreq}(b)
that $\abs{s_0(N)}^{\ell-1}\sim\lambda^N$.
Putting the obtained estimates together, we get
$\abs{s_0(N)}\alpha_{s_0(N)}\sim\eps^{-1/\ell}$,
which provides the polynomial factor in part~(a).
The estimate obtained is valid for the dominant coefficient
of the Melnikov function. To get the analogous estimate for the splitting
function, one has to bound the corresponding coefficient of the
error term in~(\ref{eq:Melniapprox}), showing that it is also exponentially
small and dominated by the main term in the approximation.
This works as in \cite{DelshamsG04}, where the case of the golden number
was considered, and we omit the details here.

The proof of part~(b) can be carried out in similar terms.
For the second dominant harmonic, we get an exponentially small estimate
with the function $h_2(\eps)$, defined in~(\ref{eq:h2}).
This estimate is also valid if one considers the whole sum in~(b),
since the terms of this sum can be bounded by a geometric series,
and hence it can be estimated by its dominant term
(see \cite{DelshamsG04} for more details).
\qed

%----------
%----------

\subsection{Study of the functions $h_1(\eps)$ and $h_2(\eps)$}\label{sect:h1}

To conclude the proof of Theorem~\ref{thm:main}, we show in this section
the different properties of the function $h_1(\eps)$
for the quadratic and cubic cases,
and establish a condition, fulfilled
in all cases~(\ref{eq:24numbers}--\ref{eq:cubicgolden}),
ensuring that $h_2(\eps)\le h_1(\eps)$
for any $\eps$.

\begin{lemma}
\ \btm
\item[\rm(a)]
In the quadratic case $\ell=2$, the function $h_1(\eps)$
is $4\ln\lambda$-periodic
in $\ln\varepsilon$, with $\min h_1 (\varepsilon) =A_0$ and
$\max h_1(\varepsilon)=A_1$, with constants
\beq\label{eq:quad_A0A1}
  A_0=1,
  \qquad
  A_1=\frac12\p{\frac1{\sqrt\lambda}+\sqrt\lambda}.
\eeq
\item[\rm(b)]
In the cubic case $\ell=3$, the function $h_1(\eps)$ satisfies the bound
$0< A_0^- \leq h_1 (\varepsilon) \leq A_1^+$,
with constants
\beq\label{eq:cub_A0A1}
  A_0^-=(1-\delta)^{1/3},
  \quad
  A_1^+=\frac{(1+\delta)^{1/3}}3
    \pq{2\p{\frac{\sqrt\lambda+1}{2\lambda}}^{1/6}
    +\p{\frac{2\lambda}{\sqrt\lambda+1}}^{1/3}}.
\eeq
\etm
\end{lemma}

\proof
We use that the functions $g^*_n(\eps)$ and the values $\eps^*_n$ satisfy
the following scaling properties:
\beq\label{eq:scaling}
  g^*_{n+1}(\eps)=\p{\frac{b_{n+1}}{b_n}}^{\,1/\ell}
    \cdot g^*_n\p{\frac{\eps^*_n}{\eps^*_{n+1}}\cdot\eps},
  \qquad
  \eps^*_{n+1}
  =\p{\frac{b_n}{b_{n+1}\lambda^\ell}}^{2/(\ell-1)}\cdot\eps^*_n\,.
\eeq

In the quadratic case ($\ell=2$), we have $b_n=b_{n+1}=1$,
and hence~(\ref{eq:scaling}) becomes
\[
  g^*_{n+1}(\eps)=g^*_n(\lambda^4\,\eps),
  \qquad
  \eps^*_{n+1}=\frac{\eps^*_n}{\lambda^4}\,,
\]
where we have $(\eps^*_n)$ as a geometric sequence,
and the functions $g^*_n(\eps)$ are just translations of the initial one,
if we use a logarithmic scale for $\eps$ (see Figure~\ref{fig:qn-1-2}).
It is easy to check that the intersection between the graphs of the
functions~$g^*_n$ and $g^*_{n+1}$ takes place at $\eps'_n:=\eps^*_n/\lambda^2$.
Thus, for $\eps\in\pq{\eps'_n,\eps'_{n-1}}$
we have $N(\eps)=n$, and hence $h_1(\eps)=g^*_n(\eps)$.
We can obtain $h_1(\eps)$ from any given interval $\pq{\eps'_n,\eps'_{n-1}}$
by extending it as a $4\ln\lambda$-periodic function
of $\ln\eps$, and it is clear that its minimum and maximum values are
$h_1(\eps^*_n)=1$ and $h_1(\eps'_n)=h_1(\eps'_{n-1})=A_1$ respectively.

The cubic case ($\ell=3$) becomes more cumbersome, because the function
$h_1(\eps)$ is not periodic in $\ln\eps$, due to the oscillating
quantities $b_n$, $b_{n+1}$ in~(\ref{eq:scaling});
notice that
$\ds\frac{1-\delta}{1+\delta}<\frac{b_{n+1}}{b_n}<\frac{1+\delta}{1-\delta}$.
Nevertheless, we are going to obtain a periodic upper bound
for the function $h_1(\eps)$.
Let us introduce the functions
\[
  g^+_n(\eps):=G(\eps;\eps^+_n,1+\delta),
  \qquad
  \varepsilon^+_n:= \frac{D_0}{K_{j_0}^{\,3}}
  \cdot\frac{1}{(1+\delta)\lambda^{3n}},
\]
obtained from~(\ref{eq:gn}), by replacing the oscillatory factors $b_n$
by the constant $1+\delta$ (and taking $\ell=3$).
One can check that the graphs of the functions $g^*_n(\eps)$ and $g^+_n(\eps)$
have no intersection if $b_n<1+\delta$, and coincide if $b_n=1+\delta$,
which implies that $g^*_n(\eps)\le g^+_n(\eps)$ for any $\eps$.
As in~(\ref{eq:h1}), we can define
\[
  h^+_1(\eps):=\min_{n\ge0}g^+_n(\eps),
\]
and it is clear that we have the upper bound
$h_1(\eps)\le h^+_1(\eps)$ for any $\eps$.
By a similar argument to the one used in the quadratic case,
we can establish the periodicity in $\ln\eps$ of the function $h^+_1(\eps)$,
and we find its maximum value.
Indeed, one can check from the expression~(\ref{eq:defG}), with $\ell=3$,
that the graphs of the functions $g^+_n$ and $g^+_{n+1}$
intersect at
\[
  \eps'_n:=\p{\frac{\sqrt\lambda+1}{2\lambda}}^2\cdot\eps^+_n,
  \qquad
  \textrm{with \ $g^+_n(\eps'_n)=g^+_{n+1}(\eps'_n)=A_1^+$.}
\]
Thus, we have $h^+_1(\eps)=g^+_n(\eps)$ for $\eps\in\pq{\eps'_n,\eps'_{n-1}}$,
and we can extend it as a $3\ln\lambda$-periodic function in $\ln\eps$,
whose maximum value is $A_1^+$,
which provides an upper bound for $h_1(\eps)$.
On the other hand, since $g^*_n(\eps)\ge b_n^{\,1/3}$ it is clear that
$h_1(\eps)\ge A^-_0$.
\qed

\bremark
In general, for the cubic case the function $h_1(\eps)$ is not periodic
in $\ln\eps$, but we can conjecture that it is quasiperiodic due to the
oscillating quantities $b_n$ introduced in~(\ref{eq:bn}).
In fact, using the approximation~(\ref{eq:phi}) for the angle $\phi$,
we have $b_{n+22}\simeq b_n\,$, which gives $66\ln\lambda$
as an approximate period for $h_1(\eps)$ (see Figure~\ref{fig:cub_h1}).
\eremark

As said in Section~\ref{sect:dominant},
the function $h_1(\eps)$ is related to the dominant harmonic
among the primary resonances, corresponding to the integer vector $s_0(N)$,
with $N=N(\eps)$.
In order to ensure that this harmonic provides the maximal splitting distance,
we need that $h_1(\eps)\le g_k(\eps)$ also for secondary harmonics $k$.
To have this inequality,
the separation between the primary and secondary resonances
has to be large enough.
Recalling that the separations $B_0$, $B_0^-$
for quadratic and cubic frequencies
were defined in~(\ref{eq:quad_B0}) and~(\ref{eq:cub_B0}) respectively,
we impose the \emph{``separation condition''}:
\beq\label{eq:B0A1}
  \begin{array}{ll}
    B_0\ge A_1     &\textrm{in the quadratic case;}
  \\[3pt]
    B_0^-\ge A_1^+ &\textrm{in the cubic case.}
  \end{array}
\eeq
A numerical exploration of this condition, among quadratic frequency vectors
given by purely periodic continued fractions,
indicates that the 24~cases considered in~(\ref{eq:24numbers})
are all the ones satisfying it.
On the other hand, the condition is fulfilled in the case of the
cubic golden vector considered in~(\ref{eq:cubicgolden}),
since $B_0^-\approx1.1824$ and $A_1^+\approx1.0909$,
as can be checked from the numerical data given
in Section~\ref{sect:cubicgolden}.

Notice that, under the separation condition~(\ref{eq:B0A1}),
we have $h_2(\eps)<h_1(\eps)$, unless $\eps$ is very close to
some concrete values where there is a change in the dominant harmonic given
by $N=N(\eps)$
(in the quadratic case, this happens for $\eps$ close
to the geometric sequence $(\eps'_n)$; see Figure~\ref{fig:qn-1-2}).

Now we can complete the proof of Theorem~\ref{thm:main}.
Indeed, according to Proposition~\ref{prop:Mdom}
the separation condition implies that for any $\eps$ sufficiently small,
we have the estimate
\[
  \max\limits_{\theta\in \mathbb{T}^\ell} |\M(\theta)|\sim\abs{\M_{s_0(N)}},
\]
since the coefficient of the dominant harmonic $k=s_0(N)$, $N=N(\eps)$,
is greater or equal than the sum of all other harmonics.

\small
%\bibliographystyle{alpha}
%\bibliography{splqc}

\def\noopsort#1{}

\end{document}